\documentclass[11pt]{amsart}
\usepackage{latexsym,xspace,enumerate,amsfonts,amsmath,amssymb,amscd}
\usepackage{xypic,color,tikz,hyperref,syntonly,MnSymbol,slashed,paralist,stmaryrd}
\usepackage[capitalize]{cleveref}
\usepackage[arrow, matrix, curve]{xy} 

\newcommand{\sfl}{\mathrm{sf}}
\newcommand{\calO}{\mathcal O}
\newcommand{\e}{\epsilon}
\newcommand{\app}{\mathrm{app}}
\newcommand{\HH}{\mathbb H}

\newcommand{\Z}{\mathbb{Z}}

\newcommand{\R}{\mathbb{R}}
\newcommand{\C}{\mathbb{C}}

\newcommand{\GL}{\mr{GL}}
\newcommand{\SL}{\mr{SL}}

\newcommand{\U}{\mr{U}}
\newcommand{\SU}{\mr{SU}}

\newcommand{\PSL}{\mr{PSL}}

\newcommand{\mf}{\mathfrak}
\newcommand{\mr}{\mathrm}
\newcommand{\mb}{\mathbf}
\newcommand{\mc}{\mathcal}

\newcommand{\End}{\mathop{\rm End}\nolimits}

\newcommand{\Id}{\mathop{\rm Id}\nolimits}

\renewcommand{\ker}{\mathop{\rm ker}\nolimits}

\renewcommand{\Re}{\mathop{\rm Re}\nolimits}
\renewcommand{\Im}{\mathop{\rm Im}\nolimits}

\newcommand{\Tr}{\mathop{\rm Tr}\nolimits}

\newcommand{\del}{\partial}
\newcommand{\delb}{\bar\partial}

\newcommand{\loc}{\mr{loc}\,}
\newcommand{\note}[1]{\marginpar{\raggedright\if@twoside\ifodd\c@page\raggedleft\fi\fi\sf\scriptsize \red{RMK: #1}}}

\newcommand\red[1]{\textcolor{red}{#1}}

\newcommand{\be}{\begin{equation}}
\newcommand{\ben}{\begin{equation}\nonumber}
\newcommand{\ee}{\end{equation}}
\newcommand{\bp}{\begin{para}}
\newcommand{\ep}{\end{para}}
\newcommand{\bps}{\begin{paras}}
\newcommand{\eps}{\end{paras}}
\newcommand{\calC}{\mathcal C}

\newtheorem{proposition}{\textbf{Proposition}}
\newtheorem{lemma}[proposition]{\textbf{Lemma}}

\newtheorem{theorem}[proposition]{\textbf{Theorem}}

\theoremstyle{definition}
\newtheorem{definition}{\textbf{Definition}}

\newtheorem*{example*}{\textbf{Example}}
\newtheorem*{remark}{\textbf{Remark}}
\newtheorem*{theorem*}{\textbf{Theorem}}

\newcounter{para}[section]
\newenvironment{para}[2][]{\refstepcounter{para}\noindent\ignorespaces{\bf #1\thepara. #2.} \rmfamily}{\noindent\ignorespacesafterend\bigskip}
\newenvironment{paras}[1]{\noindent\ignorespaces{\bf #1.} \rmfamily}{\noindent\ignorespacesafterend\bigskip}

\numberwithin{proposition}{section}
\numberwithin{definition}{section}
\newcommand{\fid}{\mr{fid}}
\newcommand{\irr}{\mr{irr}}
\begin{document}
\title{Limiting configurations for solutions of Hitchin's equation}
\date{\today}

\author{Rafe Mazzeo}
\address{Department of Mathematics, Stanford University, Stanford, CA 94305 USA}
\email{mazzeo@math.stanford.edu}

\author{Jan Swoboda}
\address{Mathematisches Institut der LMU M\"unchen\\Theresienstra{\ss}e 39\\D--80333 M\"unchen\\ Germany}
\email{swoboda@math.lmu.de}

\author{Hartmut Wei\ss{}} 
\address{Mathematisches Seminar der Universit\"at Kiel\\ Ludewig-Meyn Stra{\ss}e 4\\ D--24098 Kiel\\ Germany}
\email{weiss@math.uni-kiel.de}

\author{Frederik Witt} 
\address{Mathematisches Institut der Universit\"at M\"unster\\ Einsteinstra{\ss}e 62\\ D--48149 M\"unster\\ Germany}
\email{frederik.witt@uni-muenster.de}

\thanks{RM supported by NSF Grant DMS-1105050, JS supported by DFG Grant Sw 161/1-1}

\begin{abstract}
We review recent work on the compactification of the moduli space of Hitchin's self-duality equation. We study the degeneration 
behavior near the ends of this moduli space in a set of generic directions by showing how limiting configurations
can be desingularized. Following ideas of Hitchin, we can relate the top boundary stratum of this space of limiting configurations
to a Prym variety. A key r\^ole is played by the family of rotationally symmetric solutions to the self-duality equation on $\C$, 
which we discuss in detail here. 
\end{abstract}

\maketitle

\tableofcontents
%
%
\section{Introduction}
The moduli space of Higgs bundles, introduced by Hitchin~\cite{hi87} and Simpson~\cite{si88}, is a well investigated object in algebraic geometry 
and topology. We wish here to study it from the viewpoint of Riemannian geometry. Hitchin showed that there exists a natural hyperk\"ahler metric 
on the smooth locus of the moduli space; in many cases the moduli space has no singularities and the metric is complete. Its asymptotics are
still not well understood, and we survey here some recent approaches to a set of questions about the behavior of this metric near the ends of this
moduli space. 

\medskip

There are several reasons to study this metric carefully.  The first is to understand the $L^2$-cohomology of this space. Hausel proved \cite{ha99} that
the image of the compactly supported cohomology in the ordinary cohomology vanishes, leading him to conjecture that the $L^2$-cohomology of 
the Higgs bundle moduli space must vanish. This was made in analogy with Sen's conjecture about the $L^2$-cohomology of the monopole moduli 
spaces~\cite{se94}.  Hitchin proved a rather general result \cite{hi00} showing that under conditions satisfied in both these cases, the $L^2$-cohomology
vanishes outside the middle degree. Hausel's conjecture remains open.  Following, for example, the approach of \cite{hhm05}, an understanding
of this middle-degree cohomology relies on some finer knowledge of the metric structure at infinity. 

\medskip

However, this is part of a much broader picture concerning hyperk\"ahler metrics on algebraic completely integrable systems. Indeed, 
the work of Gaiotto, Moore and Neitzke~\cite{gmn10,gmn13} hints at an asymptotic development of this hyperk\"ahler metric $g$, where 
the leading term is a so-called {\em semiflat} metric and the correction terms decay at increasingly
fast exponential rates. The exponents and coefficients of these correction terms are described in terms of expressions coming from a wall-crossing 
formalism, but these are unfortunately a priori divergent.   Clarifying this circle of ideas is a high priority. 

\medskip

Our goal here is to review the main result of~\cite{msww14}, which constructs a dense open subset near infinity in the moduli space of 
Hitchin's self-duality equations. The degeneration behavior of generic solutions is captured by the notion of {\em limiting configurations}.
These constitute a family of singular solutions to the self-duality equations which give a geometric realization of the elements of the
top stratum in the compactification of the moduli space. As a second step, we present a desingularization theorem for limiting configurations.
We present here an alternate description of these limiting configurations, different than the one one given in \cite{msww14}; the approach
here was communicated to us by Hitchin~\cite{hi14}, and we are grateful to him for allowing us to use it here.  We review this desingularization 
result and conclude with a sketch of what these results
indicate about the asymptotic behavior of the hyperk\"ahler metric; complete proofs of this will appear in a subsequent paper. 

\medskip

The most pressing question is to understand the metric asymptotics in all the remaining ``non-generic'' cases. This parallels the story 
of the moduli space of $\SU(2)$-monopole metrics of charge $k$ studied by Bielawski. Our genericity condition corresponds to a multi-monopole 
breaking up into $k$ monopoles of charge one, cf.\ ~\cite{bi95}. In this ``free'' region the natural hyperk\"ahler metric is asymptotic to the 
so-called {\em Gibbons-Manton metric}~\cite{bi98}. For the general case which involves so-called clusters (where monopoles break up into 
monopoles of smaller charge but not necessarily equal to one), Bielawski found hyperk\"ahler metrics approximating the natural hyperk\"ahler 
metric at an exponential rate~\cite{bi08}.

\bigskip

\centerline{\textbf{Acknowledgements}}

\smallskip

This project was initiated and conducted in the research SQuaRE ``Nonlinear analysis and special geometric structures'' at the American Institute of 
Mathematics in Palo Alto. It is a pleasure to thank AIM for their hospitality and support. JS gratefully acknowledges the kind hospitality of the 
Department of Mathematics at Stanford University during a several months research visit in 2013. HW would like to thank the organizers of 
the S\'eminaire de Th\'eorie Spectrale et G\'eom\'etrie for the invitation and the opportunity to speak. Finally, the authors would like to thank Olivier
Biquard, Vicente Cort\'es and Nigel Hitchin for useful discussions. 
%
%
\section{Holomorphic bundles with Higgs fields}
In this section we review some relevant background. A more complete introduction can be found, for example, 
in the appendix in~\cite{wgp08}.  For generalities on hermitian holomorphic vector bundles see ~\cite{ko87}. 
%
\subsection{Stable bundles}\label{stab.bund}
Let $X$ be a closed Riemann surface, i.e.\ a compact (orientable) surface endowed with a complex structure. To streamline 
the discussion we always assume that the genus $\gamma$ of $X$ is at least $2$. We also fix a complex vector bundle 
$E\to X$ of rank $r=r(E)$ and degree $d=d(E)$. The pair $(r,d)$ determines $E$ as a smooth bundle. Given the holomorphic 
structure on $X$ it is natural to classify the holomorphic structures on $E$. These are characterized in terms of 
{\em pseudo-connections}, i.e., $\C$-linear maps $\delb_E:\Omega^0(E)\to\Omega^{0,1}(E)$ such that $\delb_E(fs)=
\delb f\otimes s+f\delb_Es$ for any complex-valued smooth function $f$ and section $s$ of $E$. Given such a pseudo-connection we can define a complex horizontal subspace of the tangent space of $E$. Together with the complex structure of the fibres we therefore obtain an almost complex structure on $TE$. The integrability condition for this is $\delb_E\circ\delb_E=0$ (see for instance \cite[Proposition I.3.7]{ko87}) which on a Riemann surface holds trivially for dimensional reasons. Note in passing that any pseudo-connection arises as the $(0,1)$-part of a full covariant derivative $d_A=\del_A+\delb_A$. A complex 
automorphism $g\in\mc G^c:=\mr{Aut}(E)$ acts on a holomorphic structure by
$$
\delb_E^g:=g^{-1}\circ\delb_E\circ g.
$$
The associated moduli space
$$
\mc N_{r,d}=\{\mbox{pseudo-connections }\delb_E\}/\mc G^c
$$
is not Hausdorff in general. To obtain a well-behaved quotient one must restrict attention to the subclass of {\em stable} bundles, 
following ideas from Mumford's geometric invariant theory~\cite{mu65,ne78}. For a general complex vector bundle, we define
the {\em slope} of $E$ to be the number $\mu(E)=d(E)/r(E)$. We say that $E$ is {\em (slope-)stable} if and only if 
$$
\mu(F)<\mu(E).
$$
for any proper non-trivial holomorphic subbundle $F$ of $E$. The {\em moduli space of stable bundles} $\mc N^s_{r,d}$ is 
a smooth quasi-projective variety of (complex) dimension $1+r^2(\gamma-1)$; this space is projective if $r$ and $d$ are coprime. 

\medskip

A completely different approach to the moduli space of stable bundles was developed by Atiyah, Bott and 
Donaldson~\cite{ab83,do83}. This builds upon the seminal work of Narasimhan and Seshadri~\cite{ns65}, 
who proved that a stable holomorphic structure on $E$ is equivalent to a {\em projectively flat unitary connection}. 
More concretely, fix a hermitian metric $H$ on $E$. A unitary connection $A$ is projectively flat if the induced connection 
on the associated prinicipal $\mr{PU}(r)$-bundle is flat. This may also be described directly in terms of the curvature 
$F_A$ of $A$ by the condition
\begin{equation}\label{gau.sta.bun}
F_A=-i\mu\cdot\omega.
\end{equation}
Letting $\mc G=\{g\in\mc G^c\mid gg^*=\Id_E\}$ denote the space of {\em unitary} gauge transformations, then 
there is a corresponding moduli space 
\begin{equation}\label{iso.sta.bun}
\mb N^{\irr}_{r,d}=\{A\mbox{ irreducible solution to }\eqref{gau.sta.bun}\}/\mc G.
\end{equation}
When $(r,d)$ are coprime again, this is a smooth manifold of (real) dimension $2+2r^2(\gamma-1)$ and is diffeomorphic to $\mc N^s_{r,d}$.

\begin{remark}\label{dep.her.met}
If $(E,\delb_E,H)$ is a hermitian holomorphic vector bundle, there is a unique connection $A=A(H,\delb_E)$, called the {\em Chern connection},
on $E$ which is unitary and satisfies $\delb_A=\delb_E$. If the underlying holomorphic structure is clear we shall simply write $A_H$ or even $A$ for the Chern connection associated with $H$. This gives an action of $\mc G^c$ on $H$-unitary connections by 
setting $A^g=A(H,\delb_E^g)$ (see \cite[Chapter VII]{ko87}). Donaldson showed that if $(E,\delb)$ is stable, then we can find 
in any complex gauge orbit of the Chern connection a projectively flat unitary connection. This connection is not in general 
the same as the Chern connection for $H$. It is unique only up to $H$-unitary gauge transformations. The precise choice of 
the hermitian background metric $H$ is immaterial. Indeed, complex gauge transformations act transitively on the space of 
hermitian metrics, so a change in the metric will only affect the representative of the orbit solving \cref{gau.sta.bun}, but not 
the moduli space itself. Equivalently, we can regard the holomorphic structure $\delb_E$ as being fixed and then consider 
\cref{gau.sta.bun} as an equation on $H$ via the assignement $H\mapsto A(H,\delb)$. Thus, starting with an arbitrary metric 
$H$, we find a complex gauge transformation $g$ -- unique up to $H$-unitary transformations -- such that the Chern 
connection $A(H^g,\delb_E)$ solves \cref{gau.sta.bun}.
\end{remark}
%
\subsection{Higgs bundles}\label{higgs.bundles}
The cotangent bundle of a complex manifold carries a natural holomorphic symplectic structure, i.e., a closed non-degenerate 
holomorphic $2$-form. If this is a real-analytic K\"ahler manifold, then there exists a hyperk\"ahler metric in a neighbourhood 
of the zero section of the cotangent bundle which generates this holomorphic symplectic structure~\cite{fe01,vk99} 
(see also \cref{hyp.met} for the definition of hyperk\"ahler metrics). It is straightforward, using the linearization of
the equation $\delb_E\circ\delb_E=0$ to compute that $T_E\mc N_{r,d}^s=H^1(X,\End(E))$. By Serre duality, the fibres of 
this cotangent bundle are 
$$
T^*_E\mc N_{r,d}=H^1(X,\End(E))^*\cong H^0(X,\End(E)\otimes K),
$$
where $K\to X$ is the canonical line bundle of $X$. Sections of $T_E^*\mc N_{r,d}$, i.e., holomorphic bundle maps 
$\Phi:E\to E\otimes K$, are called {\em Higgs fields}.  A {\em Higgs bundle} is a pair consisting of a holomorphic 
vector bundle $(E,\delb_E)$ and a Higgs field $\Phi$ on it. 

\medskip

As for stable bundles, this picture has a gauge theoretic interpretation given by Hitchin's {\em self-duality equations}. We
describe these now. Fix a hermitian metric $H$ on $E$. We consider pairs $(A,\Phi)$, where $A$ is a unitary connection 
and $\Phi\in\Omega^{1,0}(\End(E))$ is an (a priori smooth) Higgs field. The equations we require these to satisfy are 
\begin{equation}\label{self.dual.equ}
F_A+[\Phi\wedge\Phi^*]=-i\mu\cdot\omega,\quad\delb_A\Phi=0.
\end{equation}
These arose as a dimensional reduction of the instanton equation on $\R^4$. The term ``Higgs field'' was coined in analogy 
with the three-dimensional counterpart of \cref{self.dual.equ}, the so-called {\em Bogomolny equations}. Solutions are 
absolute minimizers of a dimensionally reduced Euclidean Yang-Mills-Higgs functional in dimension 
three~\cite{jt80}. 

\medskip

The second equation of~\eqref{self.dual.equ} states that $\Phi$ is holomorphic with respect to the holomorphic structure defined 
by $\delb_A$. Hence a solution of \eqref{self.dual.equ} specifies a Higgs bundle $(E,\delb_A,\Phi)$. Omitting the Higgs field, 
we simply recover \cref{gau.sta.bun}. Conversely, any Higgs bundle $(E,\Phi)$ where $E$ is stable in the sense of \cref{stab.bund} 
arises as an irreducible solution of \cref{self.dual.equ}. Here irreducibility means that $(E,A,\Phi)$ does not admit a decomposition 
into a direct sum of two hermitian bundles $(E_1\oplus E_2,A_1\oplus A_2,\Phi_1\oplus\Phi_2)$ with unitary connections $A_i$ 
and Higgs fields $\Phi_i\in\Omega^{1,0}(\End(E_i))$.   On the other hand, not every Higgs bundle which arises from an irreducible 
solution of \cref{self.dual.equ} necessarily has an underlying stable holomorphic vector bundle. To capture {\it all} the irreducible 
solutions $(A,\Phi)$ we therefore generalize the stability condition. Thus we say that $(E,\Phi)$ is a {\em stable Higgs bundle} if 
for any proper non-trivial holomorphic subbundle $F$ of $E$ with $\Phi(F)\subset F\otimes K$ we have $\mu(F)<\mu(E)$. 
This reduces to the usual notion of stability when $\Phi=0$. The moduli space
$$
\mc M^s_{r,d}=\{(E,\Phi)\mid\mbox{ stable Higgs bundle}\}/\mc G^c
$$
is a quasi-projective variety of (complex) dimension $2+2r^2(\gamma-1)$ and contains $T^*\mc N_{r,d}$ as an open dense 
subset~\cite{hi87,ni91,si88}. In particular, for a generic Higgs bundle the underlying holomorphic vector bundle is stable 
in the sense of \cref{stab.bund}.  We refer to \cite[Section 3]{hi87} for examples of a Higgs bundle $(E,\Phi)$ where 
$E$ is not stable. On the other hand, 
$$
\mb M^{\irr}_{r,d}=\{(A,\Phi)\mbox{ an irreducible solution to }\eqref{self.dual.equ}\}/\mc G
$$
is a smooth manifold of (real) dimension $4+4r^2(\gamma-1)$ which is diffeomorphic to $\mc M^s_{r,d}$.

\begin{remark}
From another point of view we can think of the Higgs bundle moduli space as the moduli space of irreducible and projectively 
flat {\em complex} connections~\cite{do87,co88}, i.e.\
$$
\mc M^s_{r,d}=\{A\mbox{ irreducible complex solution of }\eqref{gau.sta.bun}\}/\mc G^c.
$$
We therefore obtain a Narasimhan-Seshadri type theorem for complex connections. This has been generalized from Riemann surfaces 
to higher-dimen\-sional compact K\"ahler manifolds by Simpson in his quest to parametrize the flat complex connections which 
arise in the complex variation of a Hodge structure. For higher dimensions it is necessary to impose the extra condition 
$\Phi\wedge\Phi=0$ (which is trivially satisfied on a Riemann surface).
\end{remark}

In this paper we consider a variant of the self-duality equations where we work in the {\em fixed determinant} case. If $A$ is a 
unitary connection, then its curvature $F_A$ decomposes as
$$
F_A = F_A^\perp + \frac{1}{r} \Tr(F_A) \otimes \Id_E,
$$
where $F_A^\perp \in \Omega^2(\mf{su}(E))$ is the {\em trace-free} part of the curvature and $\mf {su}(E)\subseteq\End(E)$ denotes the bundle of traceless skew-hermitian endomorphisms. The {\em central part} 
$\Tr(F_A) \in \Omega^2(i\underline{\R})$ is precisely equal to the curvature of the induced connection on $\det E$. 
Let us fix a background connection $A_0$ from now on and consider only those connections $A$ which induce the 
same connection on $\det E$ as $A_0$ does, i.e.\ $A=A_0 + \alpha$ where $\alpha \in \Omega^1(\mf{su}(E))$. In 
other words, any such $A$ is trace-free ``relative'' to $A_0$. We may now consider the pair of equations 
\be\label{hit.equ.fixed.det}
\begin{array}{rcl}
F_A^\perp+[\Phi\wedge\Phi^*] & = & 0, \\[0.4ex]
\delb_A\Phi & =& 0,
\end{array}
\ee
where $A$ is trace-free relative to $A_0$. Since the trace of a holomorphic Higgs field is constant, we may as well restrict to 
trace-free Higgs fields $\Phi \in \Omega^{1,0}(\End_0(E))$.  There always exists a unitary connection $A_0$ on $E$ such that 
$\Tr F_{A_0}=-ir(E) \omega$, and with this as background connection, a solution of \eqref{hit.equ.fixed.det} provides 
a solution to \eqref{self.dual.equ}. Of course we now need to restrict to gauge transformations of unit determinant 
$\mc G^c_0$ and $\mc G_0$ when building the moduli spaces. The precise choice of $A_0$ is immaterial in the sense that 
the moduli spaces corresponding to two such choices are isomorphic.

\medskip

In the sequel we specialize not only to the fixed determinant case, but also to rank $2$ bundles. This is the case originally 
considered by Hitchin. We denote by $\Lambda\to X$ a fixed degree $d$ line bundle; this carries a natural connection induced by $A_0$. In particular, it is holomorphic. We consider rank $2$ stable Higgs bundles $(E,\Phi)$ with $\det E=\Lambda$ (as holomorphic line bundles). The moduli space of all such bundles is denoted $\mc M_\Lambda$, where for simplicity we drop any reference to stability or irreducibility.  
%
\subsection{Parabolic Higgs bundles}
There is an extension of the definition of stability of bundles to the setting of punctured Riemann surfaces.
We recall this briefly since we make auxilliary use of this later.  {\em Parabolic bundles} were first introduced 
by Seshadri~\cite{se77},~\cite{ms80}. The corresponding notion of a parabolic Higgs bundle was introduced 
by Simpson in~\cite{si90} and developed further in~\cite{ns95,by96}.
 
\begin{definition}
A {\em parabolic structure} on a holomorphic vector bundle $E$ over $X$ with marked points $\mf p=\{p_1,\ldots,p_n\}\subset X$ 
consists of the following data: 
\begin{itemize}
	\item a filtration $E_p=F_p^1\supsetneq\ldots\supsetneq F_p^{l_p}\supsetneq0$ at $p\in\mf p$;
	\item a system of associated {\em weights} $0\leq\alpha_1<\ldots<\alpha_{l_p}<1$.
\end{itemize}
The {\em parabolic degree} of $E$ is defined as 
$$
\mr{pardeg}(E)=d(E)+\sum_{p_i\in\mf p}\alpha_i(\dim F_i-\dim F_{i+1}).
$$
A {\em parabolic Higgs bundle} $(E,\Phi)$ consists of a holomorphic bundle $E$ with parabolic structure and a 
Higgs field which is nilpotent with respect to the flags at the marked points $p_i\in\mf p$, i.e., 
$\Phi(F_p^i)\subset F_p^{i+1}\otimes K$.  (Alternately, one may require that $\Phi(F_j)\subset F_j\otimes K$.) 
\end{definition}

To motivate this definition, consider the punctured surface $X^\times=X\setminus\mf p$. There exists a discrete 
subgroup $\Gamma$ of $\PSL(2,\R)$ acting freely on the upper half-plane $\mathbb H^2$ so that 
$X^\times=\HH^2/\Gamma$ and $\Gamma\cong\pi_1(X^\times)$. The elements of $\Gamma$ corresponding to 
loops around the punctures are of parabolic type, which means that their extensions to $\overline{\HH^2} = 
\HH^2 \cup \partial_\infty \HH^2$ have exactly one fixed point on the boundary $\partial_\infty \HH^2$, see for 
instance~\cite[Chapter A]{bp92}. If we denote by $(\HH^2)^+$ the union of $\HH^2$ with all fixed points of parabolic 
elements in $\Gamma$, then $X=(\HH^2)^+/\Gamma$ is a compact Riemann surface with marked points $p_1,\ldots,p_n$,
where each $p_j$ is the `endpoint' of a parabolic cusp, whence the appellation parabolic bundle. The isotropy 
group of a parabolic cusp under the action of $\Gamma$ is cyclic. Let $a_i$ denote the loop in $\pi_1(X^\times)=\Gamma$ 
induced by a generator of the cyclic group corresponding to $p_i$. Then a unitary representation $\rho:\Gamma\to\U(r)$ 
gives rise to matrices
\begin{equation}\label{matrix.para.cusp}
\rho(a_i)=\left(\begin{array}{ccc}\exp(2\pi i\alpha_1)&&0\\&\ddots&\\0&&\exp(2\pi i\alpha_r)\end{array}\right)
\end{equation}
for a suitably chosen unitary basis $e_1,\ldots, e_r$ of $E$ near $p_i$. This gives the weighted flag. 

\begin{remark}
There exists a natural notion of stability for parabolic (Higgs) bundles, just as in the unmarked case, leading to 
moduli spaces of stable parabolic Higgs structures. 
\end{remark}
%
\subsection{Spectral curves}\label{spec.curve}
The map
$$
\det:\mc M_\Lambda\to H^0(X,K^2),\quad[E,\Phi]\mapsto\det\Phi
$$
is proper and surjective and defines the so-called {\em Hitchin fibration} of $\mc M_\Lambda$. If $q\in H^0(X,K^2)$ has simple zeroes, 
the fibre $\det^{-1}(q)$ can be identified with the Prym variety of the spectral curve associated with $q$~\cite{hi87,hi87b}.

\medskip

Since the zeroes of $q$ are assumed to be simple,
$$
\hat X_q:=\{a_x\in T^*X\mid a_x^2=q(x)\}
$$
defines a smooth embedded curve, which is called the 
{\em spectral curve} associated with $q$. We write this simply as $\hat X$ if there is no risk of confusion. The projection 
$\pi(a_x)=x$ is a twofold cover $\pi:\hat X\to X$ ramified at the divisor $q^{-1}(0)$ in $X$. By the Riemann-Hurwitz
formula, the genus of $\hat X$ is $4\gamma-3$. We denote by $\sigma:\hat X\to\hat X$ the involution which interchanges 
the sheets of $\pi$. The pull-back bundle $\pi^*K$ admits a tautological holomorphic section $a_x\in\hat X
\mapsto a_x\in K_x=(\pi^*K)_x$ which we suggestively denote by $\sqrt q$. Now let
$$
\mr{Prym}_\Lambda(\hat X)=\{\mbox{holomorphic line bundles}\  L\to\hat X\mid\det \, \pi_*L=\Lambda\otimes K^{-1}\}
$$
where $\pi_*$ denotes the {\em direct image} of the holomorphic line bundle $L$, cf.\ \cite[\S 3b]{gu67} or \cite[Chapter 2.4]{hsw99}. 
Since $\hat X\to X$ has degree two and $L \to \hat X$ is a holomorphic line bundle, $\pi_*L\to X$ is a rank $2$ holomorphic 
vector bundle. For instance, $\pi_*\mc O_{\hat X}=\mc O_X\oplus K^{-1}$ (\cite[Ex.\ IV.2.6]{ha77}). If $\Lambda=\mc O$ we recover 
the {\em Prym variety} of $\hat X$; by definition, this is the $3\gamma-3$-dimensional subvariety of the Jacobian of $\hat X$ 
defined by the property $\sigma^*L=L^*$~\cite{mu74}.

\medskip

To define Higgs bundles from elements in $\mr{Prym}_\Lambda(\hat X)$ we first choose a holomorphic square root $K^{1/2}$ of $K$. 
For $L\in\mr{Prym}_\Lambda(\hat X)$, set $M=L\otimes\pi^*K^{1/2}$ and define $E=\pi_*M$. By the {\em projection 
formula}~\cite[Lemma 10]{gu67}, multiplication by $\sqrt q$ induces a holomorphic map $M\to M\otimes\pi^*K$ which 
descends to a holomorphic map
$$
\Phi:E=\pi_*(M)\to E\otimes K=\pi_*(M\otimes\pi^*K).
$$
Since $M$ is associated with a Prym variety, this Higgs field is trace free, and satisfies $\det\Phi=q$. Note also that
$$
\det E=\det(\pi_*L)\otimes K^{-1/2}=\Lambda,
$$
whence $E$ has the right determinant. Conversely, let $(E,\Phi)$ be a Higgs bundle on $X$ with $\det E=\Lambda$ 
and $\det\Phi=q$. Pulling back the Higgs bundle to the spectral curve $\hat X_q$ we obtain a holomorphic rank $1$ 
subbundle $T \subset \ker(\pi^*\Phi-\sqrt q)$. Then $E=\pi_*M$ for $M=T\otimes\pi^*K$, see \cite[Remark 3.7]{bnr87}, 
so that $L=M\otimes\pi^*K^{-1/2}\in\mr{Prym}_\Lambda(\hat X)$. 
%
%
\section{Limiting configurations}
As the name suggests, limiting configurations are the structures which arise as limits of solutions of the Hitchin equations.
Conversely, it is possible to desingularize such limiting configurations to obtain `large' elements in the Higgs bundle moduli
space. We wish here to explain this further and discuss the existence of such limiting configurations.
%
\subsection{Motivation}
Following~\cite{hi87} we know that the function
$$
\mc M_\Lambda\to\R,\quad[A,\Phi]\mapsto\|\Phi\|^2_{L^2}
$$
is a proper Morse-Bott function. In other words, if $(A_n, \Phi_n)$ is a sequence of solutions to Hitchin's equations, and 
if the $L^2$ norms of the $\Phi_n$ are bounded, then these solutions lie in a compact subset of $\mc M_\Lambda$. 
On the other hand, fixing the Higgs field $\Phi$, then Hitchin's existence theorem guarantees the existence of a pair 
$(A_t,t\Phi_t)$ in the complex gauge orbit of $(A_H,t\Phi)$ (in particular, $\Phi_t$ is complex gauge equivalent to $\Phi$),
satisfying
$$
F_{A_t}^\perp+t^2[\Phi_t\wedge\Phi_t^*]=0,\quad\delb_{A_t}\Phi_t=0.
$$
This family of solutions approaches the end of the moduli space as $t\to\infty$.  

To get a feeling for how degenerations occur, {\em assume} that $(A_t,\Phi_t)\to(A_\infty,\Phi_\infty)$ as $t\nearrow\infty$
in $\calC^\infty_{\loc}$ (note that we have normalized the Higgs field $t\Phi_t$ in the solution by dividing by $t$).
Then it necessarily holds that
$$
[\Phi_\infty\wedge\Phi_\infty^*]=0,\quad\delb_{A_\infty}\Phi_\infty=0.
$$
In particular, $\Phi_\infty$ is normal, and hence unitarily diagonalizable. Therefore, at a point $p\in X$ where $\det\Phi_\infty\in H^0(X,K^2)$ 
vanishes, this order of vanishing must be at least two; order one vanishing occurs if and only if $\Phi$ is nilpotent at $p$. However, for 
generic Higgs fields, $\det\Phi$ has {\em simple zeroes} (we also say that $\Phi$ is {\em simple}). This means that $\Phi_\infty$
must be singular at these zeroes. Depending on the point of view one takes, cf.\ \cref{dep.her.met}, either the hermitian metric 
degenerates (so that normality is no longer defined at $p$) or else the holomorphic structure breaks down (so $\det\Phi_\infty$ is 
not holomorphic on all of $X$). Either way, it seems reasonable to expect that $F_{A_t}$ concentrates near $(\det\Phi_\infty)^{-1}(0)$ 
as $t\nearrow\infty$. This is consistent with recent results by Taubes~\cite{taa,tab}, who investigated the analogous degeneration
behavior for a three-dimensional analogue of the self-duality equation. 
%
\subsection{The fiducial solution}
We next study a class of solutions to Hitchin's equation on $\C$ which are rotationally symmetric in an appropriate sense. 
We learned about this family of fiducial solutions from Andy Neitzke, who in joint work with Gaiotto and Moore \cite{gmn13} 
described its basic properties. We are grateful to him for explaining this to us in detail. The symmetry reduces Hitchin's 
equation to a Painlev\'e type III ordinary differential equation, and from that perspective this solution can be traced back 
to work of Mason and Woodhouse \cite{mawo93}.  

\medskip

Our paper \cite{msww14} presents a derivation of these solutions and their basic properties, but we do so again here
from a different and simpler perspective. Indeed, as we explain here, the fiducial solutions are obtained rather easily if one 
interprets Hitchin's equations as equations for the hermitian metric when restricted to the complex gauge orbit of a 
configuration $(A,\Phi)$ with $\bar\partial_A\Phi=0$. In general, fixing a background hermitian metric $H_0$ on the vector bundle $E$, 
we may identify an arbitrary hermitian metric with a hermitian endomorphism field $H$; in the fixed determinant case, this field
also satisfies $\det H = 1$. Let $\ast_H$ denote the adjoint taken with respect to the hermitian metric $H$. Then 
\begin{equation}
F^\perp_{A} + \bar\partial_{A}(H^{-1} \partial_{A} H) + [\Phi \wedge \Phi^{*_H}]=0
\label{hmh}
\end{equation}
if and only if $(A_H,\Phi)^g$ solves Hitchin's equation, where $g$ is a complex gauge transformation which is 
determined by solving $H^{-1}=gg^*$. Note that $g$ is uniquely determined up to right multiplication by 
a unitary gauge transformation.

\medskip

Fix the holomorphic normal form for a rank-$2$ Higgs bundle over $\C$ whose determinant has a simple zero in $z=0$. 
More precisely, consider the holomorphically trivial rank-$2$ vector bundle over $\C$ and the holomorphic Higgs field 
\[
\Phi =  \varphi \, dz, \quad  \varphi(z)=\begin{pmatrix} 0 & 1 \\ z & 0 \end{pmatrix};
\]
this has determinant $-zdz^2$. Note that if $H_0$ is the standard constant hermitian metric, then the Chern connection is the trivial flat connection,
denoted here by $A$.

\medskip

We ask first if there exists a configuration $(A_\infty,\Phi_\infty)$ in the complex gauge orbit of $(A, \Phi)$ satisfying 
the purely algebraic equation
\begin{equation}\label{algebraic}
[\Phi_\infty \wedge \Phi_\infty^*] =0 \Longleftrightarrow [\Phi \wedge \Phi^{*_{H_\infty}}]=0; 
\end{equation}
here $H_\infty$ is a hermitian metric and $\ast_{H_\infty}$ is the adjoint with respect to $H_\infty$. 

\medskip

We say that a hermitian metric is rotationally symmetric if its representation $H$ relative to the fixed metric 
$h_0$ is rotationally symmetric, i.e., 
\[
H(r)= \begin{pmatrix}
 \alpha(r) & b(r) \\ \bar b(r) & \beta(r)
 \end{pmatrix},
\]
where $\alpha, \beta$ are real valued, with $\alpha>0$ and $\alpha\beta-|b|^2=1$. A straightforward calculation shows that 
\begin{equation}
[\Phi \wedge \Phi^{*_{H}}] = \begin{pmatrix} \alpha^2 - |z|^2\beta^2 + 2i \Im b^2 z & 2 \alpha b -2   \beta \bar b \bar z\\
 -2  \alpha b z + 2  |z|^2 \beta \bar b& |z|^2\beta^2 - \alpha^2 - 2i \Im b^2 z
\end{pmatrix} dz \wedge d\bar z. 
\end{equation}
We wish to choose $H = H_\infty$ so that this vanishes. Since $\alpha, \beta \neq 0$, setting all entries here equal to $0$ 
implies that $b \equiv 0$, and then that $\alpha  = r^{1/2}$, $\beta = r^{-1/2}$. In other words, 
\[
[\Phi \wedge \Phi^{*_{H_\infty}}]=0 \Longleftrightarrow H_\infty = \begin{pmatrix}
 r^{1/2} & 0 \\ 0 & r^{-1/2}
 \end{pmatrix}
\]
This solution is singular at $z=0$. Recalling that the curvature of any Hermitian metric $H$ is given by
\begin{equation}
\bar \partial (H^{-1} \partial H)= -\begin{pmatrix} \partial_{\bar z} (\beta \partial_z \alpha - b \partial_z \bar b) & \partial_{\bar z} (\beta \partial_z b - b \partial_z \beta) \\ \partial_{\bar z} (\alpha \partial_z \bar b-\bar b \partial_z \alpha) & \partial_{\bar z} (\alpha \partial_z \beta-\bar b \partial_z b) \end{pmatrix} d z \wedge d \bar z,
\end{equation}
then in particular,
\[
\bar \partial (H_{\infty}^{-1} \partial H_\infty)= - \begin{pmatrix}
 \partial_{\bar z} (r^{-1/2} \partial_z r^{1/2}) & 0\\
 0 & \partial_{\bar z} (r^{1/2} \partial_z r^{-1/2})
 \end{pmatrix} d z \wedge d \bar z = 0.
\]

This proves the following surprising fact:

\begin{lemma}
The unique hermitian metric $H_\infty$ satisfying $[\Phi \wedge \Phi^{*_{H_\infty}}]=0$ on $\C^\times$ is flat and therefore solves the decoupled version of Hitchin's equation
\[
\bar \partial (H_\infty^{-1}\partial H_\infty) = 0, \quad [\Phi \wedge \Phi^{*_{H_\infty}}]=0.
\]
\end{lemma}

Setting 
\[
g_\infty= \begin{pmatrix}
 r^{-1/4} & 0\\
 0 & r^{1/4}
 \end{pmatrix},
\]
then one has $H_\infty^{-1}=g_\infty g_\infty^*$, so from the fact that $(A_\infty^{\fid},\Phi_\infty^{\fid}) = (A,\Phi)^{g_\infty}$,
we obtain the expression for the limiting fiducial solution 
\[
A^\fid_\infty=\frac{1}{8}\begin{pmatrix}1&0\\0&-1\end{pmatrix}\left(\frac{dz}{z}-\frac{d\bar z}{\bar z}\right),\quad\Phi^\fid_\infty=
\begin{pmatrix}0&r^{1/2}\\ z r^{-1/2}&0\end{pmatrix}dz. 
\]

Now let us look for nonsingular rotationally symmetric solutions $H_t$ of Hitchin's equation
\begin{eqnarray}\label{eq:HitchinPDE}
\bar \partial (H_t^{-1} \partial H_t) + t^2[\Phi \wedge \Phi^{*_{H_t}}]=0, \quad 0 < t < \infty.
\end{eqnarray}
We see that $[\Phi \wedge \Phi^{*_{H_t}}]$ is rotationally symmetric if and only if 
\[
H_t = \begin{pmatrix}\alpha_t&0\\0&\alpha_t^{-1}\end{pmatrix}
\]
for some function $\alpha_t=\alpha_t(r) > 0$. With this reduction, Hitchin's equation reduces to the ODE
\begin{eqnarray}\label{eq:HitchinODE}
\partial_{\bar z}(\alpha_t^{-1}\partial_z\alpha_t) - t^2(\alpha_t^2-|z|^2\alpha_t^{-2})=0.
\end{eqnarray}
We let
\begin{eqnarray*}
\alpha_t=e^{h_t + \tfrac 12 \log(r)}
\end{eqnarray*}
where $h_t=h_t(r)$ is real-valued. Then \eqref{eq:HitchinODE} is equivalent to 
\begin{multline*}
0=\partial_{\bar z}\partial_z(h_t+\tfrac{1}{2}\log(r))-t^2(re^{2h_t}-re^{-2h_t})\\
=\tfrac{1}{4}\Bigl(\frac{d^2}{dr^2}+\frac{1}{r}\frac{d}{dr}\Bigr)(h_t+\tfrac{1}{2}\log(r))-t^2(re^{2h_t}-re^{-2h_t})\\
=\tfrac{1}{4}\big(h_t''+r^{-1}h_t'\big)-t^2(re^{2h_t}-re^{-2h_t}),
\end{multline*}
and hence finally to 
\begin{equation}\label{ODE}
h_t''+r^{-1}h_t'=8t^2r\sinh(2h_t).
\end{equation}
Now substitute $h_t(r)=\psi(\rho)$ with $\rho=\frac 83 tr^{3/2}$; this transforms \eqref{ODE} to the Painlev\'e type III equation
\begin{equation}\label{Painleve}
\psi''+\frac{\psi'}{\rho}=\frac 12 \sinh(2\psi).
\end{equation}
The properties of solutions to \eqref{Painleve} are well known, see \cite{msww14} and the references therein. There 
is a unique solution $h_t$ to \eqref{ODE} satisfying $h_t(r)+\frac{1}{2}\log(r)\to 0$ as $r\searrow 0$ and $h_t(r)\to 0$ 
as $r \nearrow \infty$.

The hermitian metric $H_t$ then equals 
\[
H_t = \begin{pmatrix} r^{1/2} e^{h_t(r)} & 0 \\ 0 & r^{-1/2} e^{-h_t(r)} \end{pmatrix},
\]
and 
\[
H_t^{-1} = g_t g_t^* \quad \mbox{for} \quad g_t= \begin{pmatrix}
 r^{-1/4}e^{-h_t(r)/2} & 0\\
 0 & r^{1/4}e^{h_t(r)/2}
 \end{pmatrix}. 
\]
This yields, finally, the entire family of desingularized fiducial solutions by the formula
$(A_t^{\fid},\Phi_t^{\fid}) = (A,\Phi)^{g_t}$, so
\[
A^\fid_t=f_t(r)\begin{pmatrix}1&0\\0&-1\end{pmatrix}\left(\frac{dz}{z}-\frac{d\bar z}{\bar z}\right),\quad\Phi_t^\fid =\begin{pmatrix} 0 & r^{1/2} e^{h_t(r)} \\ zr^{-1/2}e^{-h_t(r)} & 0 \end{pmatrix}\, dz
\]
where $f_t(r) = \frac{1}{8} + \frac{1}{4} r \partial_r h_t$.

\medskip

We have therefore arrived at the following result:

\begin{theorem}
For each $t \in (0,\infty)$ there exists a unique hermitian metric $H_t$ solving Hitchin's equation \eqref{eq:HitchinPDE} such that
\begin{enumerate}[(i)]
\item $H_t$ is rotationally symmetric, and
\item $H_t \sim H_\infty=\begin{pmatrix} r^{\frac{1}{2}}&0\\0&r^{-\frac{1}{2}}\end{pmatrix}$ as $r\nearrow \infty$.
\end{enumerate}
These metrics correspond to the fiducial solutions $(A_t^\fid,\Phi_t^\fid)$.
\end{theorem}


A byproduct of this approach is that we may now easily derive fiducial solutions corresponding to Higgs fields $\Phi$, the determinants
of which have zeroes of order greater than one.   This touches on the ongoing thesis work of Laura Fredrickson \cite{lafr}, who
is describing the behavior of families of solutions on $\C$ near this degenerate case, i.e., families where the Higgs field
is simple but which limit to these degenerate solutions. She also constructs fiducial solutions in the higher rank case.
More precisely, for $k\geq1$, consider the Higgs field
\[
\Phi = \begin{pmatrix} 0 & 1 \\ z^k & 0 \end{pmatrix} dz,
\]
which has determinant $-z^kdz^2$. Setting 
\[
H_t=\begin{pmatrix}\alpha_t&0\\0&\alpha_t^{-1}\end{pmatrix}
\]
we obtain  that
\[
[\Phi \wedge \Phi^{*_{H_t}}] = \begin{pmatrix} \alpha_t^2 - |z|^{2k}\alpha_t^{-2}   &0\\
 0 & |z|^{2k}\alpha_t^{-2} - \alpha_t^2 
 \end{pmatrix} dz \wedge d\bar z .
\]
If we let $\alpha_t=e^{h_t+\frac{k}{2}\log(r)}$, then the same calculation as above leads now to the ODE
\begin{eqnarray}\label{higher_order_ODE}
h_t''+r^{-1}h_t'=8t^2r^k\sinh(2h_t).
\end{eqnarray}
With the substitution
\begin{eqnarray*}
\rho=4tr^{\frac{k}{2}+1}/(\tfrac k2 +1)
\end{eqnarray*} 
this is equivalent to the same Painlev\'e  III equation \eqref{Painleve}
as in the case where $\det \Phi$ has a simple zero. The Higgs pair $(A_t,\Phi_t)$ corresponding to this rotationally symmetric 
hermitian metric $H_t$ is now calculated to be 
\[
A_t=f_t(r)\begin{pmatrix}1&0\\0&-1\end{pmatrix}\left(\frac{dz}{z}-\frac{d\bar z}{\bar z}\right), \quad \Phi_t=\begin{pmatrix}0&r^{k/2}e^{h_t(r)}\\z^kr^{-k/2} e^{-h_t(r)}&0\end{pmatrix}\,dz
\]
where $f_t(r) = \frac{1}{8} + \frac{1}{4} r \partial_r h_t$ and $h_t$ is a solution of \eqref{higher_order_ODE}. The limiting solution is given by
\[
A_\infty=\frac{k}{8}\begin{pmatrix}1&0\\0&-1\end{pmatrix}\left(\frac{dz}{z}-\frac{d\bar z}{\bar z}\right), \quad \Phi_\infty=\begin{pmatrix}0&r^{k/2}\\z^kr^{-k/2}&0\end{pmatrix}\,dz.
\]

The gluing result from \cite{msww14} may be implemented using this higher order fiducial solution. Details will 
appear elsewhere. 
%
\subsection{Construction of limiting configurations}\label{limi.prym}
Motivated by the preceding discussion we consider the following model for degeneration at infinity~\cite{msww14}. 

\begin{definition}\label{lim.conf.equ}
Let $(E,H,\Phi)$ be a hermitian Higgs bundle, where $\Phi$ is simple. A {\em limiting configuration} is a Higgs pair 
$(A_\infty,\Phi_\infty)$ which satisfies the decoupled self-duality equations 
\begin{equation}
F_{A_\infty}^\perp=0,\quad[\Phi_\infty\wedge\Phi_\infty^*]=0,\quad\delb_{A_\infty}\Phi_\infty=0
\end{equation}
on $X^\times$, and which agrees with $(A_\infty^\fid, \Phi_\infty^\fid)$ near each point of $(\det\Phi)^{-1}(0)$ with respect to some 
unitary frame for $E$ and holomorphic coordinate system such that $\det\Phi=-zdz^2$ near the zeroes of $\det\Phi$. 
\end{definition}

The main result for limiting configurations is this.

\begin{theorem}\label{ex.out.sol}
Let $(E,H,\Phi)$ be a Hermitian Higgs bundle with simple Higgs field. Then in the complex gauge orbit of $(A_H,\Phi)$ over 
$X\setminus(\det\Phi)^{-1}(0)$ there exists a limiting configuration $(A_\infty,\Phi_\infty)$ unique up to unitary gauge transformation. 
Up to a smooth complex gauge transformation on $X$ the limiting (singular) complex gauge transformation $g_\infty$ 
takes the form
$$
g_\infty=\left(\begin{array}{cc}|z|^{-\tfrac{1}{4}}&0\\0&|z|^{\tfrac{1}{4}}\end{array}\right)
$$
near $(\det\Phi)^{-1}(0)$. 
\end{theorem}
 
This was proved in \cite{msww14} analytically, in particular relying on the Fredholm theory of conic elliptic operators. 
We present here an alternative approach using spectral curves which was explained to us by Nigel 
Hitchin~\cite{hi14}, and we are grateful for his permission to present it here. 

\medskip

Following the notation from \cref{spec.curve}, let $\hat X$ be the spectral curve associated with $q=-\det\Phi$. Start with 
$L\in\mr{Prym}_\Lambda(\hat X)$, but now consider the holomorphic rank $2$ bundle $\hat E=L\oplus\sigma^*L$ over 
$\hat X$. Near the ramification locus (which is the zero locus of $\sqrt q$) we fix a holomorphic trivialization of $L$ 
and declare it to be unitary. This determines a hermitian metric $\hat h$ near these points, which we extend to 
all of $\hat X$ and set $\hat H=\hat h\oplus\sigma^*\hat h$. Finally, with respect to the decomposition of $\hat E$ 
we define a Higgs field by $\hat\Phi=\mr{diag}(\sqrt q,-\sqrt q)$.  This determines a hermitian Higgs bundle 
$(\hat E,\hat H,\hat\Phi)$ where $\hat\Phi$ is normal on $\hat X$ and $\hat H$ is flat near the ramification locus. 

\medskip

The $\Z_2$-action on $\hat X$ generated by $\sigma$ is covered by a representation 
$$
\tau:\Z_2\to\GL_r\C, \qquad \mbox{where}\qquad \tau(\sigma)=\left(\begin{array}{cc}-1&0\\0&1\end{array}\right)
$$
in a suitable basis $(e,e')$. Writing the action of $\sigma$ on an eigenvector as $e^{i\pi x}$, we obtain the {\em isotropy weights} 
$x=0$ and $x'=1$ for this action. Since the Higgs bundle $(\hat E,\hat H,\hat\Phi)$ is $\sigma$-invariant, it descends to $X$ 
and defines there an orbifold hermitian Higgs bundle $(E,H,\Phi)$, see~\cite{ns95} and also~\cite{bo91,fs92}. In particular, 
the underlying orbifold bundle has trivializations of the form $\hat D\times\C^2/\sigma\times\tau$ on neighbourhoods 
$D=\hat D/\sigma$ of points in $q^{-1}(0)$.  In addition, if $w$ is a coordinate on $\hat D$, then $z=w^2$ is a coordinate on $D$. 

\medskip

Orbifold Higgs bundles can be desingularized in a natural way by parabolic Higgs bundles with cusps at $\mf p=q^{-1}(0)$.
Namely, at a marked point $p$ with isotropy weights $x=0$ and $x'=1$ we consider the bundle $\mc E$ defined by
$$
E|_{X\setminus\{p\}}\cup_\Psi D\times\C^2
$$
with clutching function $\Psi$ given in local coordinates by its $\Z_2$-equivariant lifting
\begin{align*}
\hat\Psi:(\hat D\setminus\{0\})\times\C^2&\to D\times\C^2\\
(w,(v_1,v_2))&\mapsto(w^2,(w^{-x}v_1,w^{-x'}v_2))=(w^2,(v_1,w^{-1}v_2)).
\end{align*}
Thus $\mc E$ carries a natural parabolic structure which at a marked point $p\in q^{-1}(0)$ is given by 
\begin{equation}\label{flag}
\begin{array}{ccccccc}
\mc E_p&=&\C^2&\supset&\C e'&\supset&\{0\}\\
0&\leq&0&<&\frac 12&<&1
\end{array}.
\end{equation}

\medskip

To describe the holomorphic sections of $\mc E$ near $p$ we note that the invariant holomorphic sections of
$(\hat D^2\times\C^2)/(\sigma\times\tau)$ with respect to $(e,e')$ take the form $(s_1(w),s_2(w))$ with 
$s_1(w)=\tilde s_1(w^2)$ and $s_2(w)=w\tilde s_2(w^2)$. Under the clutching map $\Psi$ this gives a bounded 
holomorphic section $(\tilde s_1(z),\tilde s_2(z))$ of $\mc E$ which extends over $0$. In particular, $f=e$ and $f'=e/w$ induce a holomorphic frame for $\mc E$ on the punctured disk $D\setminus\{0\}$. It follows that 
$\pi^*\Lambda^2\mc E=\pi^*(\Lambda\otimes K^{-1})$, hence $\Lambda^2\mc E=\Lambda\otimes K^{-1}$ for 
the map $\pi^*:\mr{Jac}(X)\to\mr{Jac}(\hat X)$ between the respective Jacobians is injective \cite{mu74}. Furthermore, with respect to the frame $(f,f')$, the induced hermitian metric $\mc H$ is given by
$$
\mc H=\left(\begin{array}{cc}1&0\\0&|z|\end{array}\right).
$$
As for the Higgs field we find  
$$
\Phi(z)=\left(\begin{array}{cc}0&z\\1&0\end{array}\right)dz
$$
by \cite[Section 5A]{ns95}. Note that $\Phi$ is indeed parabolic at $z=0$. 

\medskip

To get a holomorphic vector bundle with deteminant $\Lambda$ we must as in \cref{spec.curve} twist with the 
square root $K^{1/2}$ to get $E_\infty=\mc E\otimes K^{1/2}$ with corresponding Higgs field $\Phi_\infty:E_\infty\to E_\infty\otimes K^{1/2}$. 
Restricted to $X^\times$ any line bundle is holomorphically trivial, and a hermitian metric is just a nowhere vanishing function 
$h$ with respect to this trivialization. We trivialize $K^2$ by $q$ and define $h=|q\bar q|^{-2}$. Taking $\sqrt[4]{h}$ for a trivialization of $K^{1/2}$ yields a hermitian metric which is given by $1/|z|$ near the marked points and is locally of the form $f\bar f$ for a holomorphic function $f$. In particular, this metric is flat so that the product metric $H_\infty$ on $E_\infty|_{X^\times}$ has the same curvature as $\mc H$. Near a marked point,
$$
H_\infty=\left(\begin{array}{cc}\tfrac{1}{\sqrt{|z|}}&0\\0&\sqrt{|z|}\end{array}\right)
$$
so the associated Chern connection $A_\infty$ is precisely equal to $A^\fid_\infty$. Similarly, writing $\Phi_\infty$ with respect to the 
rescaled unitary frame gives the fiducial Higgs field $\Phi^\fid_\infty$.

\medskip

To show that $(A_\infty,\Phi_\infty)$ is complex gauge equivalent to $(A_H,\Phi)$ on $X^\times$ we may assume (cf.\ \cref{dep.her.met}) that modulo a smooth complex gauge transformation over $X$ our initial hermitian metric $H$ is given by $\mr{diag}(1,1)$ with respect to the frame $(f,f')$. Then $H_\infty=g_\infty^* H g_\infty$ on $X^\times$ where 
$$
g_\infty=\left(\begin{array}{cc}|z|^{-\tfrac{1}{4}}&0\\0&|z|^{\tfrac{1}{4}}\end{array}\right)
$$
near the marked points.

\begin{remark}
This viewpoint is somehow dual to the one in~\cite{msww14}. Indeed, here we look for a limiting hermitian metric, which
then determines a holomorphic bundle $E_\infty$ endowed with a singular hermitian metric. By contrast, in~\cite{msww14} 
we look for a limiting unitary connection; this gives a smooth hermitian metric but a singular complex structure. 
\end{remark}
%
%
\section{Desingularization by gluing}
We now give a brief sketch of the main theorem in \cite{msww14}, which globalizes the 
phenomenon that the family of smooth fiducial solutions $(A_t^\fid, \Phi_t^\fid)$
desingularize the limiting fiducial solution $(A_\infty^\fid, \Phi_\infty^\fid)$.

\begin{theorem}[\cite{msww14}]\label{gluing_theorem}
Suppose that $(A_\infty, \Phi_\infty)$ is a simple limiting configuration, i.e., so that 
$\det \Phi_\infty$ has simple zeroes. Then there exists a family of solutions to the 
Hitchin equations $(A_t, t \Phi_t)$ with $A_t \to A_\infty$ and $\Phi_t \to \Phi_\infty$ 
in $\calC^\infty_\loc$ at exponential rate in $t$ on the complement of $\det \Phi_\infty^{-1}(0)$. 
\end{theorem}

These solutions are obtained by gluing $(A_\infty,\Phi_\infty)$ to the 
elements of the family $(A_t^\fid,\Phi_t^\fid)$. In fact, it is possible to prove this 
theorem either by performing the gluing as just indicated, or else by constructing a
desingularizing family of hermitian metrics $H_t$ as solutions to \eqref{hmh}.  
These two approaches are quite close, and there seems to be nothing to recommend 
one over the other, but we explain the former. 

\medskip

The steps in the proof follow a familiar pattern.  We first construct a family of approximate
solutions $A_t^\app$ and $\Phi_t^\app$ as follows.  Since the limiting configurations, by
definition, agree with the limiting fiducial solution in a neighbourhood of each marked point,
and noting that $(A_t^\fid, \Phi_t^\fid)$ converges exponentially in $t$ to $(A_\infty^\fid, \Phi_\infty^\fid)$ 
on any annulus $\{\frac12 \e \leq |z| \leq \e\}$, we can use a partition of unity to patch 
$(A_t^\fid, \Phi_t^\fid)$ to $(A_\infty, \Phi_\infty)$ to obtain $(A_t^\app, \Phi_t^\app)$. The error term, 
i.e., the deviation which measures the extent to which these do not satisfy Hitchin's equations, 
is supported in this annulus and is exponentially small in $t$.   For the second step we seek 
a small complex gauge transformation $g = \exp \gamma$ which has the property that
$(A_t^\app,t\Phi_t^\app)^g$ is an exact solution.   Note that we may assume
that $\gamma$ is a section of the bundle $ i \mf{su}(E)$ of hermitian endomorphisms of $E$
since this is transverse to the infinitesimal real gauge transformations. 

\medskip

The second step requires more work.  To set things up, fix a background connection $A_0$ and 
consider the Hitchin operator 
\[
\mc H_t(A,\Phi) = (\mc H_t^{(1)}, \mc H_t^{(2)}) = (F^\perp_A+t^2[\Phi \wedge \Phi^*],\bar\partial_A\Phi)
\]
for connections $A$ which are traceless relative to $A_0$ and traceless Higgs fields $\Phi$. Consider also the orbit map
\begin{equation}
\mathcal{O}_{(A,\Phi)}(\gamma)= (A^g,\Phi^g), \qquad g = \exp (\gamma). 
\label{orb.map}
\end{equation}
The goal is to find a pair $(A, \Phi)$ which lies in the nullspace of $\mc H_t \circ \calO_{(A,\Phi)}$. Since the condition 
$\bar{\partial}_A \Phi = 0$ is complex gauge invariant, we may disregard the second component of 
$\mc H_t$, so it suffices to find a solution of 
\begin{equation}
\mathcal F_t(\gamma) := \mathcal{H}_t^{(1)} \circ \mathcal{O}_{(A,\Phi)}(\exp (\gamma)) = 0, 
\label{eq:defmapFt}
\end{equation}
or more explicitly, 
\[
F_{A^g}^\perp+t^2[\Phi^g\wedge(\Phi^g)^*] = 0, \qquad g = \exp(\gamma).
\]

Since we may assume that $\gamma\in \Omega^0(i\mf{su}(E))$, we calculate first that
\[
\left. D\calO_{(A,\Phi)}\right|_{g=\mathrm{Id}} (\gamma) =(\bar\partial_A\gamma-\partial_A\gamma,[\Phi,\gamma])
\] 
and then that 
\[
D \mc H_t^{(1)} ( \dot A, \dot \Phi ) =  d_A (\dot A) +  t^2( [\Phi \wedge {\dot \Phi \,}^{\ast}]+[\Phi^{\ast}\wedge \dot \Phi \,])
\]
whence
\[
D\mathcal F_t (\gamma) = 
(\partial_A\bar\partial_A-\bar\partial_A\partial_A) \gamma+t^2([\Phi\wedge[\Phi,\gamma]^{\ast}]+[\Phi^{\ast}\wedge[\Phi,\gamma]]).
\]
Using the standard identities 
\[
2\bar\partial_A\partial_A=F_A-i\ast\Delta_A,\qquad2\partial_A\bar\partial_A=F_A+i\ast\Delta_A,
\]
as well as the fact that 
\[
[ \Phi \wedge [\Phi,\gamma]^*] = - [ \Phi \wedge [ \Phi^*,\gamma]],
\]
we can write 
\begin{equation}\label{eq:operatorDt}
D\mathcal F_t(\gamma) =  i\ast\Delta_A\gamma+t^2M_\Phi \gamma,
\end{equation}
where
\[
M_\Phi\gamma:=[\Phi^{\ast}\wedge[\Phi,\gamma]] - [ \Phi \wedge [ \Phi^*,\gamma]].
\]  
We study instead the operator on $\Omega^2(\mf{su}(E))$ given by 
\[
L_t(\gamma) = -i \,\ast  D\mathcal F_t (\gamma) = \Delta_A \gamma  - i \ast t^2 M_\Phi \gamma.
\]

A brief calculation, see \cite{msww14}, shows that if $\gamma \in \Omega^0(i \mf{su}(E))$, then
\begin{equation}
\langle L_t \gamma, \gamma \rangle = || d_{A} \gamma||^2 + 4 \, || [\Phi,\gamma] ||^2  \geq 0.
\label{ibp}
\end{equation}
In particular, $L_t \gamma =0$ if and only if $[\Phi,\gamma]=0$ and $d_A \gamma = 0$. 

\medskip

The main analytic part of this gluing argument is to determine the mapping properties of the inverse of $L_t$, and
to keep track of this behavior uniformly in $t$. As part of this, we must also show that $L_t$ is invertible.  To do
all of this, decompose $X$ into a union of disks $B_{\e}(p_j)$ around each $p_j$ and the complementary (or exterior) region
$X \setminus \sqcup B_{\e}(p_j)$.  Assume that the underlying metric on $X$ is Euclidean in each $B_\e(p_j)$. We then
analyze the restriction of this operator to each disk, which we denote by $L_t^{\mathrm{int}}$, with Neumann boundary 
conditions. This is done using separation of variables in polar coordinates, and is an explicit but lengthy computation. 
We see from this that $L_t^{\mathrm{int}}$ is invertible and satisfies
\[
||(L_t^{\mathrm{int}})^{-1}||_{\mathcal L(L^2, L^2)} \leq C, \qquad ||(L_t^{\mathrm{int}})^{-1}||_{\mathcal L(L^2, H^2)} \leq Ct^2. 
\]
On the other hand, $L_t^{\mathrm{ext}}$, the restriction of $L_t$ to the exterior region, is independent of $t$, so it 
suffices only to check that it has no nullspace (with Neumann boundary conditions).  Integrating by parts as in \eqref{ibp} 
shows that if $\gamma$ is in the nullspace of $L_t^{\mathrm{ext}}$ (with Neumann boundary conditons), then 
$\gamma$ is a parallel section (with respect to $A_\infty$ of the line bundle $\mathcal L_{\Phi_\infty}$ of elements 
$\gamma$ which satisfy $[\Phi_\infty, \gamma] = 0$, i.e., which commute with the limiting Higgs field $\Phi_\infty$. 
However, this line bundle is nontrivial around each $p_j$ and thus has no parallel sections, so we conclude that
$\gamma = 0$ and this nullspace is trivial.    Finally, a standard domain decomposition lemma shows that the operator 
$L_t$ is invertible on the entire surface $X$ and satisfies
\[
||L_t^{-1}||_{\mathcal L(L^2, L^2)} \leq C, \qquad ||L_t^{-1}||_{\mathcal L(L^2, H^2)} \leq Ct^2. 
\]

\medskip

The rest of the proof of the gluing theorem proceeds by expanding $\mathcal F_t$ into its first order Taylor 
approximation plus remainder:
\[
\mathcal F_t(\exp(\gamma))=\mathcal{H}_t^{(1)} (A_t^{\app},\Phi_t^{\app})+L_t\gamma+Q_t(\gamma).
\]
Since $L_t$ is invertible, the solutions of $\mathcal F_t(\exp(\gamma)) = 0$ are the same as solutions of 
\[
\gamma=-L_t^{-1}\big(\mathcal{H}_t^{(1)}(A_t^{\app},\Phi_t^{\app})+Q_t(\gamma)\big),
\]
or in other words, as fixed points of the mapping
\[
T\colon B_\rho\to H^2(i\mf{su}(E)),\quad\gamma\mapsto -L_t^{-1}\big(\mathcal{H}_t^{(1)}(A_t^{\app},\Phi_t^{\app})+Q_t(\gamma)\big),
\]
where $B_\rho$ is a ball in $H^2$. 

To show that $T$ really does have fixed points, we show that it is a contraction when $\rho$ is sufficiently small.
The key points here are that the norm of $L_t^{-1}$ grows polynomially in $t$ while the error term 
$\mathcal{H}_t^{(1)}(A_t^{\app},\Phi_t^{\app})$ decays exponentially in $t$. 

We remark that the exact same proof can be used to desingularize limiting configurations $(A_\infty, \Phi_\infty)$
where $\det \Phi_\infty$ has multiple zeroes. Indeed, we simply use the corresponding family of fiducial solutions
introduced earlier.   The more serious issue, however, is to understand families of solutions near to these ones
with multiple zeroes. As mentioned earlier, this is likely possible using the forthcoming work of Fredrickson \cite{lafr}.
%
%
\section{The hyperk\"ahler metric}\label{hyp.met}

Now let us turn to a consideration of the natural $L^2$, or Weil-Petersson, metric on $\mb M_\Lambda$, the moduli space of solutions to the 
self-duality equation~\eqref{hit.equ.fixed.det}.    As shown initially by Hitchin, this metric is of a very special type. 

\begin{definition}
A manifold $M^{4k}$ endowed with a Riemannian metric $g$ and three integrable complex structures $I_1$, $I_2$ and $I_3$ is called 
{\em hyperk\"ahler} if
\begin{itemize}
	\item $g$ is K\"ahler with respect to $I_j$, $j=1,2,3$;
	\item these complex structures satisfy the quaternion relations $I_1I_2=I_3$ etc.
\end{itemize}
\end{definition}

As Hitchin shows, $\mb M_\Lambda$ carries a natural hyperk\"ahler metric on its smooth locus; the smooth locus is the entire moduli space if $d$ is odd ~\cite{hi87}, and this metric is then complete by Uhlenbeck compactness. We recall next how this metric arises through a moment map interpretation of 
the self-duality equations. 
%
\subsection{Moment maps}
An action of a Lie group $G$ on a symplectic manifold $(M,\omega)$ is called {\em hamiltonian} if there exists a $G$-equivariant 
{\em moment map} $\mu:M\to\mf g^*$ from $M$ to the dual of the Lie algebra $\mf g$ of $G$. This map is defined by the property 
that for $v\in\mf g$, the smooth map $\mu_v(x)=\mu(x)(v)$ satisfies $d_x\mu_v=\omega(v^\sharp(x),\cdot)$, where $v^\sharp(x)$ 
denotes the induced fundamental vector field at $x$. If $G$ acts freely and $\lambda\in\mf g^*$ is central, the quotient 
$\mu^{-1}(\lambda)/G$ carries a natural symplectic structure inherited from $M$.

\medskip

This symplectic quotient construction can be adapted to hyperk\"ahler manifolds as follows. Assume that there is a hamiltonian action by $G$ 
with respect to any of the three K\"ahler forms $\omega_j=g(I_j\cdot,\cdot)$, $j=1,\,2,\,3$. We then assemble the individual symplectic moment 
maps $\mu_j$ into $\mu=(\mu_1,\mu_2,\mu_3)$, which takes values in $\R^3\otimes\mf g^*$.  It turns out that for a central element $\lambda
\in\R^3\otimes\mf g^*$, the quotient $\mu^{-1}(\lambda)/G$ inherits a natural hyperk\"ahler metric~\cite{hklr87}.

\medskip

This construction applies to our setting as follows~\cite{hi87}. The choice of a base connection $A_0$ identifies the space of unitary connections of fixed determinant
with $\Omega^1(\mf{su}(E))$. This, in turn, can be identified 
with the complex vector space $\mc A:=\Omega^{0,1}(\End_0(E))$, the complex conjugate of the space of traceless Higgs fields $\Omega^{1,0}(\End_0(E))$. 
The $L^2$-inner product
$$
\langle(\alpha,\Phi),(\alpha,\Phi)\rangle=2i\int_X\Tr(\alpha^*\alpha+\Phi\Phi^*)
$$
gives the quaternionic vector space $\mc A\times\bar{\mc A}$ the structure of a flat hyperk\"ahler manifold with K\"ahler forms $\omega_j$ 
induced by the complex structures
$$
I_1(\alpha,\Phi)=(i\alpha,i\Phi),\quad I_2(\alpha,\Phi)=(i\Phi^*,-i\alpha^*),\quad I_3(\alpha,\Phi)=(-\Phi^*,\alpha^*).
$$
Formally, the gauge group $\mc G$ acts in a hamiltonian fashion with moment maps
$$
\mu_1(A,\Phi)=F_A+[\Phi\wedge\Phi^*],\quad\mu_2(A,\Phi)=\Re(\delb_A\Phi),\quad\mu_3(A,\Phi)=\Im(\delb_A\Phi).
$$
The $L^2$-metric on $\mc A\times\bar{\mc A}$ induces a hyperk\"ahler metric on the quotient $\mu^{-1}(-\mu i\omega\otimes\Id,0,0)/\mc G
=\mb M_\Lambda$.

\begin{remark}
Note that for any $(a,b,c)\in S^2 \subset \R^3$, the complex structure $I_{(a,b,c)}:=aI_1+bI_2+cI_3$ is an isometry for $g$. However, the corresponding holomorphic structures on $\mb M_\Lambda$ are not equivalent. In fact, $(\mb M_\Lambda,\pm I_1)$ is biholomorphic to $\mc M_\Lambda$, while for any other complex structure defined by $S^2\setminus\{(\pm1,0,0)\}$, $(\mb M_\Lambda,I_{(a,b,c)})$ is biholomorphic to the moduli space of irreducible projectively flat complex connections together with its natural complex structure as discussed in the Remark of \cref{higgs.bundles}.
\end{remark}
%
\subsection{The semi-flat metric}
One of the most interesting open questions is to determine the asymptotic structure of the hyperk\"ahler metric $g$ on $\mc M_\Lambda$
which we identify with $\mb M_\Lambda$ using the complex structure $\pm I_1$. Guiding any investigation into this problem is the far-reaching 
conjectural picture by Gaiotto-Moore-Neitzke \cite{gmn10,gmn13} which describes this metric as a perturbation series off a so-called semi-flat 
metric $g_{\sfl}$. More specifically, $g_{\sfl}$ is a metric on the open dense set of the moduli space where the Higgs field $\Phi$ is simple, i.e., such 
that $\det \Phi$ has only simple zeroes. To understand the name `semiflat', we recall that $\mc M_\Lambda$ is the total space of a singular 
fibration $\det: \mc M_\Lambda \to \mathcal H \mathcal Q$, where the base is the space of holomorphic quadratic differentials on $X$; this
 is called the Hitchin fibration. The regular fibre of $\det$ is the Prym variety associated with $q$, which is a torus of real dimension 
$6\gamma - 6$, where $\gamma$ is the genus of $X$, cf.\ \cref{spec.curve}. 

Restricting to the regular part of the fibration this data constitute an algebraic completely integrable system and therefore the base 
carries a special K\"ahler structure \cite{fr99}. A special K\"ahler manifold is a K\"ahler manifold $(M,g,J,\omega)$ together with 
an additional flat torsion-free connection $\nabla$ on $TM$ satisfying $\nabla\omega=0$ and $d_\nabla J=0$, where $J$ is 
viewed as a $TM$-valued 1-form. As shown in \cite{fr99} one may now use the horizontal distribution on $T^*\!M$ to lift 
the metric $g$ to the horizontal part of a hyperk\"ahler metric $g_{\sfl}$ on $T^*\!M$. The restriction of $g_{\sfl}$ to each
cotangent fibre is then a flat metric, which explains the terminology. Moreover, since locally the total space of the integrable system 
is represented as the quotient of $T^*\!M$ by a family of lattices (parallel with respect to $\nabla$), the semiflat metric descends
to a metric on the regular part of the Hitchin fibration, again denoted by $g_{\sfl}$. This metric is incomplete, however, and needs 
to be corrected in order to extend over the singular fibres. A description of this correction process using a wall-crosing formalism 
is one of the main achievements of the work of Gaiotto, Moore and Neitzke. 


The first application of our new construction of `large' solutions in $\mc M_\Lambda$ shows that the hyperk\"ahler metric
$g$ is indeed well approximated by the semi-flat metric $g_{\sfl}$ far out in the ends of the moduli space. The way we seek to 
approach this is as follows: As shown in \cite{msww14}, the space of limiting configurations associated with $q \in \mathcal H 
\mathcal Q$ with simple zeroes is a torus of real dimension $6\gamma - 6$, just like the Prym variety is. A family of limiting 
configurations $(A_\infty(s),\Phi_\infty(s))$ associated with a curve of holomorphic quadratic differentials $q(s)$ may now, 
using Theorem \ref{gluing_theorem}, be perturbed into a family of solutions $(A(s),t\Phi(s))$ of Hitchin's equation for sufficiently
large $t$. After a suitable gauge choice, the derivative of this family with respect to $s$ is a tangent vector to $\mc M_\Lambda$ 
which is vertical with respect to the Hitchin fibration. The $L^2$-metric may now be directly evaluated on this tangent vector,
and is easily seen to be an exponentially small (in $t$) correction to the vertical part of the semiflat metric as $t \to \infty$. 
Horizontal tangent vectors may be obtained by varying $q \in \mc H \mc Q$, leading likewise to a comparison between the 
horizontal parts of the metrics.

%
%
%

\end{document}